\documentclass[12pt,a4paper]{amsart}
\usepackage{amssymb,amsfonts,amsthm,amsmath,graphicx,xcolor,csquotes}
\usepackage[scr]{rsfso}
\usepackage{chngcntr}
\usepackage{enumitem}

\theoremstyle{plain}

\numberwithin{theorem}{section}
\numberwithin{equation}{section}
\numberwithin{statement}{section}
\numberwithin{lemma}{section}
\numberwithin{definition}{section}
\numberwithin{conjecture}{section}
\numberwithin{corollary}{section}
\begin{document}
\title[VPV Identities, $\mathbf{x^y = y^x}$ and $\mathbf{x^y y^x = v^w w^v}$]{VPV Identities related to $\mathbf{x^y = y^x}$ and $\mathbf{x^y y^x = v^w w^v}$}
\author{Geoffrey B Campbell}
\address{Mathematical Sciences Institute,
         The Australian National University,
         Canberra, ACT, 0200, Australia}

\email{Geoffrey.Campbell@anu.edu.au}



\begin{abstract}
We cover rational and integer solutions for the equations $\mathbf{x^y = y^x}$ and $\mathbf{x^y y^x = v^w w^v}$. The former equation solutions go back to Euler, and the latter equation solutions appear to be new. Another definitely new related topic is application of VPV identities to give transforms of infinite products from our solutions. The present paper is essentially chapter 28 of the author's book appearing in June 2024.
\end{abstract}

\maketitle

\section{VPV identities in square hyperpyramid regions} \label{S:Intro VPV hyperpyramids}

In this note, which is essentially chapter 28 of the book (see Campbell \cite{gC2024}) we give some surprising Visible Point Vector identities related to rational and integer solutions of $x^y = y^x$ and also what seem to be new solutions for the equation $x^y y^x = v^w w^v$ and associated Visible Point Vector infinite product transforms. In the 18th century Christian Goldbach and later Leonhard Euler \cite{lE1748}\index{Euler, L.} famously solved the equation
\begin{equation} \label{xy.01}
  x^y = y^x
\end{equation} \label{xy.02}
in rational numbers. All possible solutions are well known to be
\begin{equation}
  x = \left(1 + \frac{1}{n}\right)^n, \quad
  y = \left(1 + \frac{1}{n}\right)^{n+1},
\end{equation}
for all positive integers $n$. The best reference for this may be the 1990 paper by M\'{a}rta Sv\'{e}d \cite{Sved1990}.

\section{Solutions to \MakeLowercase{$\mathbf{x^y y^x = v^w w^v}$} in rationals and integers}

In this section, we show there is an infinite set of nontrivial solutions to
\begin{equation} \label{xy.03}
  x^y \; y^x = v^w \; w^v
\end{equation}
in rational numbers $x$, $y$, $v$, $w$. This was recently given in Campbell\index{Campbell, G.B.} \cite{gC2022}.  However, defining what is meant by a trivial solution may be open to discussion. For a clearly nontrivial example, we see that $x = \frac{1}{3}, \; y = \frac{1}{6}, \; v = \frac{1}{2}, \; w = \frac{4}{3}$, gives us
\begin{equation} \nonumber
  \left(\frac{1}{3} \right)^{\frac{1}{6}} \; \left(\frac{1}{6} \right)^{\frac{1}{3}}
  = \left(\frac{1}{2} \right)^{\frac{4}{3}} \; \left(\frac{4}{3} \right)^{\frac{1}{2}} .
\end{equation}

Another evidently nontrivial example $x = \frac{1}{2}, \; y = \frac{1}{3}, \; v = \frac{1}{2}, \; w = \frac{4}{3}$, gives the solution

\begin{equation} \nonumber
  \left(\frac{1}{2} \right)^{\frac{1}{3}} \; \left(\frac{1}{3} \right)^{\frac{1}{2}}
  = \left(\frac{1}{2} \right)^{\frac{4}{3}} \; \left(\frac{4}{3} \right)^{\frac{1}{2}}.
\end{equation}

Noting the right sides of the above two examples are the same, we get a third non-trivial case $x = \frac{1}{2}, \; y = \frac{1}{3}, \; v = \frac{1}{3}, \; w = \frac{1}{6}$, with the solution
\begin{equation} \nonumber
  \left(\frac{1}{2} \right)^{\frac{1}{3}} \; \left(\frac{1}{3} \right)^{\frac{1}{2}}
  = \left(\frac{1}{3} \right)^{\frac{1}{6}} \; \left(\frac{1}{6} \right)^{\frac{1}{3}}.
\end{equation}

Note that trivial solutions exist if $x = v = 1$. Another solution that may be trivial (or not) is $x = \frac{1}{2}, \; y = \frac{1}{2}, \; v = \frac{1}{2}, \; w = 1$. Another possibly trivial solution is $x = 4, \; y = \frac{3}{2}, \; v = 1, \; w = \frac{81}{2}$.
A definitely trivial set of solutions is $x = a^b \; b^a, \; y = 1, \; v = a, \; w = b$, for specific fixed positive integers $a$, $b$, and $c$.

For rationals $a$, $b$ and $c$, an infinite set of rational solutions can be found setting $y=ax, \; v=bx, \; w=cx$ for any rational $x$ and substituting into $x^y y^x = v^w w^v$. This gives the algebraic solution: For $a+1\neq b+c$ with $a \neq 0$,
\begin{equation} \label{xy.04}
  x = \left(\frac{b^c c^b}{a}\right)^{\frac{1}{a-b-c+1}}, \;
  y = a\left(\frac{b^c c^b}{a}\right)^{\frac{1}{a-b-c+1}}, \;
  v = b\left(\frac{b^c c^b}{a}\right)^{\frac{1}{a-b-c+1}}, \;
  w = c\left(\frac{b^c c^b}{a}\right)^{\frac{1}{a-b-c+1}}.
\end{equation}
We seek rational solutions, so it's clear that choosing $a=b+c$ will solve our quest. Hence, an infinite number of rational solutions to $x^y y^x = v^w w^v$ are given by integer values of $b$ and $c$ with $b+c\neq0$, such that
\begin{equation} \label{xy.05}
  x = \left(\frac{b^c c^b}{b+c}\right), \;
  y = b^c c^b, \;
  v = b\left(\frac{b^c c^b}{b+c}\right), \;
  w = c\left(\frac{b^c c^b}{b+c}\right).
\end{equation}
Some verifiable random cases of this are shown here.

(b,c;x,y,v,w):
(1,	1;	$\frac{1}{2}$,	1,	$\frac{1}{2}$,	$\frac{1}{2}$),
(2,	1;	$\frac{2}{3}$,	2,	$\frac{4}{3}$,	$\frac{2}{3}$),
(3,	1;	$\frac{3}{4}$,	3,	$\frac{9}{4}$,	$\frac{3}{4}$),

(4,	1;	$\frac{4}{5}$,	4,	$\frac{16}{5}$,	$\frac{4}{5}$),
(1,	2;	$\frac{2}{3}$,	2,	$\frac{2}{3}$,	$\frac{4}{3}$),
(2,	2;	4,	16,	8,	8),
(3,	2;	$\frac{72}{5}$,	72,	$\frac{216}{5}$,	$\frac{144}{5}$),

(6,	2;	288,	2304,	1728,	576),
(1,	3;	$\frac{3}{4}$,	3,	$\frac{3}{4}$,	$\frac{9}{4}$),
(3,	3;	$\frac{3^5}{2}$,	$3^6$,	$\frac{3^6}{2}$,	$\frac{3^6}{2}$),

(6,	3;	17496,	157464,	104976,	52488),
(2,	4;	$\frac{2^7}{3}$,	$2^8$,	$\frac{2^8}{3}$,	$\frac{2^9}{3}$),

(5,	3;	$\frac{30375}{8}$,	30375,	$\frac{151875}{8}$,	$\frac{91125}{8}$), etc.

\bigskip

WolframAlpha tells us that each side of
\begin{equation*}
  288^{2304} \times 2304^{288} = 1728^{576} \times 576^{1728}
\end{equation*}
is approximately $6.843 \times 10^{6634}$ which is a 6635 digit number in base 10.

\bigskip

It also tells us that each side of
\begin{equation} \nonumber
  17496^{157464} \times 157464^{17496} = 104976^{52488} \times 52488^{104976}
\end{equation}

is approximately $9.8661 \times 10^{759039}$ which is a 759040 digit number in base 10.

\bigskip

Next, we remark that the results for $x^y = y^x$ are seen as a classical curiosity that never went beyond recreational mathematics. However, $x^y y^x = v^w w^v$ seems to have not been studied nor considered in the extended literature and the set of rational solutions given here do not seem to have been published anywhere, other than by Campbell\index{Campbell, G.B.} in 2022 \cite{gC2022}. Also, an infinite number of positive integer solutions to $x^y y^x = v^w w^v$ evidently exist, but this is not established yet. There are probably an entire class of functional equations\index{Functional equation} similar to those considered in this chapter, susceptible to rational or integer solutions, such as perhaps $v^x x^z z^v = u^v w^x y^z$, each having their own possible transforms as cases of VPV identities in combination.

\section{VPV Identity transforms using \MakeLowercase{$\mathbf{x^y = y^x}$}}

The Visible Point Vector 2D identities for the first quadrant lattice points we recall are for $|X|<1$ and $|Y|<1$ and nonzero denominators,

  \begin{equation}   \label{xy.06}
    \prod_{\substack{ (j,k)=1 \\ j,k \geq 1}} \left( \frac{1}{1-X^j Y^k} \right)^{\frac{1}{k}}
    = \left(  \frac{1}{1-Y}  \right)^{\frac{1}{1-X}},
  \end{equation}

  \begin{equation}   \label{xy.07}
    \prod_{\substack{ (j,k)=1 \\ j,k \geq 1}} \left(1-X^j Y^k \right)^{\frac{1}{k}}
    = \left(  1-Y \right)^{\frac{1}{1-X}}
  \end{equation}

These can be applied to the previous section's rational and integer solutions, by transforming each solution of $x^y = y^x$ and $x^y y^x = v^w w^v$.

We start with the first six solutions of $x^y = y^x$ in the form of $x = \left(1 + \frac{1}{n}\right)^n$, $y = \left(1 + \frac{1}{n}\right)^{n+1}$, listed as follows,

\begin{eqnarray}
\nonumber  n=1 &\Rightarrow& x = 2^1, \; y = 2^{2}; \\
\nonumber  n=2 &\Rightarrow& x = \left(\frac{3}{2}\right)^2, \; y = \left(\frac{3}{2}\right)^{3}; \\
\nonumber  n=3 &\Rightarrow& x = \left(\frac{4}{3}\right)^3, \; y = \left(\frac{4}{3}\right)^{4}; \\
\nonumber  n=4 &\Rightarrow& x = \left(\frac{5}{4}\right)^4, \; y = \left(\frac{5}{4}\right)^{5}; \\
\nonumber  n=5 &\Rightarrow& x = \left(\frac{6}{5}\right)^5, \; y = \left(\frac{6}{5}\right)^{6}; \\
\nonumber  n=6 &\Rightarrow& x = \left(\frac{7}{6}\right)^6, \; y = \left(\frac{7}{6}\right)^{7}.
\end{eqnarray}

Hence, for example, the $n=1$ only integer nontrivial solution,
\begin{equation}   \label{xy.08}
  \left( 2 \right)^{\left(4\right)} = \left( 4 \right)^{\left(2\right)}
\end{equation}
can be written in the right side of (\ref{xy.06}) form
\begin{equation}   \label{xy.09}
   \left(\frac{1}{1-\frac{1}{2}}\right)^{\left(\frac{1}{1-\frac{3}{4}}\right)}
  = \left(\frac{1}{1-\frac{3}{4}}\right)^{\left(\frac{1}{1-\frac{1}{2}}\right)},
\end{equation}
where $X= \frac{1}{2}$, $Y= \frac{3}{4}$.

Applying (\ref{xy.07}) to (\ref{xy.09}) gives us the transform

\begin{equation}   \label{xy.10}
  \prod_{\substack{ (j,k)=1 \\ j,k \geq 1}} \left\{1-\left(\frac{1}{2}\right)^j \left(\frac{1}{4}\right)^k \right\}^{\frac{1}{k}}
    =
  \prod_{\substack{ (j,k)=1 \\ j,k \geq 1}} \left\{1-\left(\frac{1}{4}\right)^j \left(\frac{1}{2}\right)^k \right\}^{\frac{1}{k}}.
  \end{equation}

Likewise, the $n=2$ case leads to the solution
\begin{equation}   \label{xy.11}
  \left( \left(\frac{3}{2}\right)^2  \right)^{\left(\frac{3}{2}\right)^{3}} = \left( \left(\frac{3}{2}\right)^{3}  \right)^{\left(\frac{3}{2}\right)^2}.
\end{equation}
This solution can then be written in the right side of (\ref{xy.06}) form
\begin{equation}   \label{xy.12}
   \left(\frac{1}{1-Y}\right)^{\left(\frac{1}{1-X}\right)}
  = \left(\frac{1}{1-X}\right)^{\left(\frac{1}{1-Y}\right)},
\end{equation}
where $X= 1-\left(\frac{2}{3}\right)^2$, $Y= 1-\left(\frac{2}{3}\right)^3$.

Applying (\ref{xy.07}) to (\ref{xy.11}) gives us the transform

\begin{equation}   \label{xy.13}
  \prod_{\substack{ (j,k)=1 \\ j,k \geq 1}} \left\{1-\left(1-\left(\frac{2}{3}\right)^2\right)^j \left(1-\left(\frac{2}{3}\right)^3\right)^k \right\}^{\frac{1}{k}}
   \end{equation}
 \begin{equation} \nonumber
    =
  \prod_{\substack{ (j,k)=1 \\ j,k \geq 1}} \left\{1-\left(1-\left(\frac{2}{3}\right)^3\right)^j \left(1-\left(\frac{2}{3}\right)^2\right)^k \right\}^{\frac{1}{k}}
  \end{equation}
  which is the same as
\begin{equation} \nonumber
    \left( \left(\frac{2}{3}\right)^2  \right)^{\left(\frac{3}{2}\right)^{3}} = \left( \left(\frac{2}{3}\right)^{3}  \right)^{\left(\frac{3}{2}\right)^2}.
\end{equation}

Again, we apply this same approach to the general rational solutions $x = \left(1 + \frac{1}{n}\right)^n$, $y = \left(1 + \frac{1}{n}\right)^{n+1}$
\begin{equation}   \label{xy.14}
  \left( \left(1 + \frac{1}{n}\right)^n  \right)^{\left(1 + \frac{1}{n}\right)^{n+1}} = \left( \left(1 + \frac{1}{n}\right)^{n+1} \right)^{\left(1 + \frac{1}{n}\right)^n}.
\end{equation}
This solution can then be written in the right side of (\ref{xy.06}) form
\begin{equation}   \label{xy.15}
   \left(\frac{1}{1-Y}\right)^{\left(\frac{1}{1-X}\right)}
  = \left(\frac{1}{1-X}\right)^{\left(\frac{1}{1-Y}\right)},
\end{equation}
where $X= 1-\left(\frac{n}{n+1}\right)^{n}$, $Y= 1-\left(\frac{n}{n+1}\right)^{(n+1)}$.

Applying (\ref{xy.07}) to (\ref{xy.11}) gives us the transform valid for all positive integers $n$,

\begin{equation}   \label{xy.16}
  \prod_{\substack{ (j,k)=1 \\ j,k \geq 1}} \left\{1-\left(1-\left(\frac{n}{n+1}\right)^n\right)^j \left(1-\left(\frac{n}{n+1}\right)^{(n+1)}\right)^k \right\}^{\frac{1}{k}}
  \end{equation}

\begin{equation} \nonumber
  \quad  =
  \prod_{\substack{ (j,k)=1 \\ j,k \geq 1}} \left\{1-\left(1-\left(\frac{n}{n+1}\right)^{(n+1)}\right)^j \left(1-\left(\frac{n}{n+1}\right)^n\right)^k \right\}^{\frac{1}{k}}
  \end{equation}
  which is the same as
\begin{equation} \nonumber
    \left( \left(\frac{n}{n+1}\right)^n  \right)^{\left(\frac{n+1}{n}\right)^{(n+1)}}
   = \left( \left(\frac{n}{n+1}\right)^{(n+1)}  \right)^{\left(\frac{n+1}{n}\right)^n}.
\end{equation}

\section{VPV Identity transforms using \MakeLowercase{$\mathbf{x^y y^x = v^w w^v}$}}

Recall earlier we stated an algebraic solution of $x^y y^x = v^w w^v$:

For $a+1\neq b+c$ with $a \neq 0$,
\begin{equation} \label{xy.17}
  x = \left(\frac{b^c c^b}{a}\right)^{\frac{1}{a-b-c+1}}, \;
  y = a\left(\frac{b^c c^b}{a}\right)^{\frac{1}{a-b-c+1}}, \;
  v = b\left(\frac{b^c c^b}{a}\right)^{\frac{1}{a-b-c+1}}, \;
  w = c\left(\frac{b^c c^b}{a}\right)^{\frac{1}{a-b-c+1}}.
\end{equation}

We then selected $a=b+c$ to find an infinite rational set of solutions $x$, $y$, $v$, and $w$, which appeared incidentally to yield an infinite set of positive integer solutions too.

Following a similar process to the previous section, we can deduce that
\begin{equation}   \label{xy.18}
   \left(\frac{1}{1-X}\right)^{\left(\frac{1}{1-Y}\right)} \left(\frac{1}{1-Y}\right)^{\left(\frac{1}{1-X}\right)}
  = \left(\frac{1}{1-V}\right)^{\left(\frac{1}{1-W}\right)} \left(\frac{1}{1-W}\right)^{\left(\frac{1}{1-V}\right)},
\end{equation}
where for non-zero denominators
\begin{equation} \label{xy.19}
  X= \frac{x-1}{x}, Y= \frac{ax-1}{ax}, V= \frac{bx-1}{bx}, W= \frac{cx-1}{cx},
  \quad with \quad x = \left(\frac{b^c c^b}{a}\right)^{\frac{1}{a-b-c+1}}.
\end{equation}

Taking the reciprocal of both sides of (\ref{xy.18}) we then apply (\ref{xy.07}) to arrive for non-zero denominators at
\begin{equation}   \label{xy.20}
  \prod_{\substack{ (j,k)=1 \\ j,k \geq 1}}
  \left\{\left[1-\left(\frac{x-1}{x}\right)^j \left(\frac{ax-1}{ax}\right)^k\right]
         \left[1-\left(\frac{ax-1}{ax}\right)^j \left(\frac{x-1}{x}\right)^k\right]\right\}^{\frac{1}{k}}
  \end{equation}
\begin{equation} \nonumber
  \quad  =
  \prod_{\substack{ (j,k)=1 \\ j,k \geq 1}}
  \left\{\left[1-\left(\frac{bx-1}{bx}\right)^j \left(\frac{cx-1}{cx}\right)^k\right]
         \left[1-\left(\frac{cx-1}{cx}\right)^j \left(\frac{bx-1}{bx}\right)^k\right]\right\}^{\frac{1}{k}},
  \end{equation}
\begin{equation} \nonumber
  \textmd{in which} \; x = \left(\frac{b^c c^b}{a}\right)^{\frac{1}{a-b-c+1}}, \; \textmd{provided that} \; a+1\neq b+c \; \textmd{and} \; a\neq0.
  \end{equation}

As an example, consider
\begin{equation} \nonumber
  288^{2304} \times 2304^{288} = 1728^{576} \times 576^{1728},
\end{equation}
which is the same as
\begin{equation} \nonumber
  \left( \frac{1}{1-\frac{287}{288}} \right)^{\left( \frac{1}{1-\frac{2303}{2304}} \right)} \times \left( \frac{1}{1-\frac{2303}{2304}} \right)^{\left( \frac{1}{1-\frac{287}{288}} \right)}
     \end{equation}
 \begin{equation} \nonumber
= \left( \frac{1}{1-\frac{1727}{1728}} \right)^{\left( \frac{1}{1-\frac{575}{576}} \right)} \times \left( \frac{1}{1-\frac{575}{576}} \right)^{\left( \frac{1}{1-\frac{1727}{1728}} \right)};
\end{equation}
and then consider the example
\begin{equation} \nonumber
  17496^{157464} \times 157464^{17496} = 104976^{52488} \times 52488^{104976},
\end{equation}
which is the same as
\begin{equation} \nonumber
  \left( \frac{1}{1-\frac{17495}{17496}} \right)^{\left( \frac{1}{1-\frac{157463}{157464}} \right)} \times \left( \frac{1}{1-\frac{157463}{157464}} \right)^{\left( \frac{1}{1-\frac{17495}{17496}} \right)}
  \end{equation}
 \begin{equation} \nonumber
= \left( \frac{1}{1-\frac{104975}{104976}} \right)^{\left( \frac{1}{1-\frac{52487}{52488}} \right)} \times \left( \frac{1}{1-\frac{52487}{52488}} \right)^{\left( \frac{1}{1-\frac{104975}{104976}} \right)}.
\end{equation}
Applying (\ref{xy.20}) then, we get the two VPV transform identities as follows.

\begin{equation}   \label{xy.21}
  \prod_{\substack{ (j,k)=1 \\ j,k \geq 1}}
  \left\{\left[1-\left(\frac{287}{288}\right)^j \left(\frac{2303}{2304}\right)^k\right]
         \left[1-\left(\frac{2303}{2304}\right)^j \left(\frac{287}{288}\right)^k\right]\right\}^{\frac{1}{k}}
  \end{equation}
\begin{equation} \nonumber
  \quad  =
  \prod_{\substack{ (j,k)=1 \\ j,k \geq 1}}
  \left\{\left[1-\left(\frac{1727}{1728}\right)^j \left(\frac{575}{576}\right)^k\right]
         \left[1-\left(\frac{575}{576}\right)^j \left(\frac{1727}{1728}\right)^k\right]\right\}^{\frac{1}{k}};
  \end{equation}
and
\begin{equation}   \label{xy.22}
  \prod_{\substack{ (j,k)=1 \\ j,k \geq 1}}
  \left\{\left[1-\left(\frac{17495}{17496}\right)^j \left(\frac{157463}{157464}\right)^k\right]
         \left[1-\left(\frac{157463}{157464}\right)^j \left(\frac{17495}{17496}\right)^k\right]\right\}^{\frac{1}{k}}
  \end{equation}
\begin{equation} \nonumber
  \quad  =
  \prod_{\substack{ (j,k)=1 \\ j,k \geq 1}}
  \left\{\left[1-\left(\frac{104975}{104976}\right)^j \left(\frac{52487}{52488}\right)^k\right]
         \left[1-\left(\frac{52487}{52488}\right)^j \left(\frac{104975}{104976}\right)^k\right]\right\}^{\frac{1}{k}}.
  \end{equation}

\end{document}